\documentclass[11pt,a4paper,twoside]{scrartcl}
\usepackage[english]{babel}

\usepackage{amsfonts}
\usepackage{amsmath}
\usepackage{amssymb}
\usepackage{amsthm}

\usepackage{graphicx}
\graphicspath{{pics/}}

\usepackage{enumerate}
\usepackage{verbatim}

\usepackage{
	float,		
	subcaption, 
	xfrac,		
	colortbl,
	mathtools,  
	pifont,  	
	mathrsfs	
}

\usepackage[symbol]{footmisc}
\usepackage{enumitem} 

\usepackage{hyperref}
\usepackage{todonotes}

\usepackage{pgfplots}
\pgfplotsset{compat=1.8}
\usepackage{tikz}
\usepackage{tikz-cd}

\newtheorem{theorem}{Theorem}[section]

\newtheorem{remark}[theorem]{Remark}

\newtheorem{example}[theorem]{Example}
\newtheorem{corollary}[theorem]{Corollary}

\newtheorem{alg}[theorem]{Algorithm}

\newenvironment{Remark}{\goodbreak \begin{remark}\slshape}{\end{remark}}
\newenvironment{Example}{\goodbreak \begin{example}\slshape}{\end{example}}

\makeatletter
\def\imod#1{\allowbreak\mkern10mu({\operator@font mod}\,\,#1)}
\makeatother

\makeatletter 
\let\c@algorithm\c@theorem  
\makeatother

\newenvironment{algorithm}[1]{\goodbreak~\begin{alg}[#1]~\vspace{-9pt}~\\
		\rule{\linewidth}{0.5pt}~\\}{\vspace{-9pt}~\\
		\rule{\linewidth}{0.5pt}~\end{alg}}

\numberwithin{equation}{section}
\numberwithin{table}{section}
\numberwithin{figure}{section}

\newcommand{\e}{\mathrm e}
\renewcommand{\i}{\mathrm i}
\newcommand{\sinc}{\mathrm{sinc}}
\renewcommand{\b}{\boldsymbol}

\newcommand{\R}{\mathbb R}
\newcommand{\C}{\mathbb C}
\newcommand{\Z}{\mathbb Z}
\newcommand{\N}{\mathbb N}
\newcommand{\T}{\mathbb T}

\newcommand{\I}{\mathcal I}

\long\def\symbolfootnote[#1]#2{\begingroup%
	\def\thefootnote{\fnsymbol{footnote}}\footnote[#1]{#2}\endgroup}

\allowdisplaybreaks

\title{Nonequispaced fast Fourier transforms for bandlimited functions}

\date{}
\author{Melanie Kircheis\footnotemark[1] \and Daniel Potts\footnotemark[3]}

\hypersetup{pdfauthor={Melanie Kircheis, Daniel Potts},
	pdftitle={Nonequispaced fast Fourier transforms for bandlimited functions},
	pdfsubject={},
	pdfcreator = {pdflatex},
	plainpages=false,
	pdfstartview=FitH,       
	pdfview=FitH,            
	pdfpagemode=UseOutlines, 
	bookmarksnumbered=true, 
	bookmarksopen=false,     
	bookmarksopenlevel=0,   
	colorlinks=true,       
	linkcolor=black,
	citecolor=black,
	urlcolor=black}

\begin{document}
	
\maketitle

\begin{abstract}
	In this paper we consider the problem of approximating function evaluations~\mbox{$f(\b x_j)$} at given nonequispaced points~$\b x_j$, \mbox{$j=1,\dots N$}, of a band\-limited function from given values~\mbox{$\hat{f}(\b k)$}, \mbox{$\b k\in \I_{\b M}$}, of its Fourier transform.
	Note that if a trigonometric polynomial is given, it is already known that this problem can be solved by means of the nonequispaced fast Fourier transform (NFFT).
	In other words, we introduce a new NFFT-like procedure for band\-limited functions, which is based on regularized Shannon sampling formulas.
	
	\smallskip
	
	\noindent AMS \emph{Subject Classifications}: \text{
		65Txx, 
		65T50, 
		94A20 
	}
\end{abstract}

\footnotetext[1]{Corresponding author: melanie.kircheis@math.tu-chemnitz.de, Chemnitz University of
	Technology, Faculty of Mathematics, D--09107 Chemnitz, Germany}
\footnotetext[3]{potts@math.tu-chemnitz.de, Chemnitz University of
	Technology, Faculty of Mathematics, D--09107 Chemnitz, Germany}


\section{Introduction}

The nonequispaced fast Fourier transform (NFFT) is a fast algorithm to evaluate a trigonometric polynomial
\begin{align*}
	f(\b x) = \sum_{\b k \in \I_{\b M}} \hat{f}_{\b k}\, \e^{2\pi\i \b k \b x},
	\quad \b x\in\T^d,
\end{align*}
with given Fourier coefficients \mbox{$\hat f_{\b k}\in\C$}, \mbox{$\b k\in\I_{\b{M}}$},
at given non\-equi\-spaced points \mbox{$\b x_j\in\T^d$}, \mbox{$j=1,\dots,N$}, \mbox{$N\in\N$}, 
where for 
\mbox{$M\in 2\N$} we define the index set~\mbox{$\I_{\b M} \coloneqq \Z^d \cap \left[-\tfrac{M}{2},\tfrac{M}{2}\right)^d$} 
with cardinality \mbox{$|\I_{\b M}| = M^d$}, and~\mbox{$\T^d \coloneqq \R^d / \Z^d$}, \mbox{$d\in\N$}, denotes the $d$-dimensional torus.

In this paper we focus on the analogous problem for band\-limited functions, where we aim to approximate evaluations~\mbox{$f(\b x_j)$}, \mbox{$j=1,\dots,N$}, of a function
\begin{align}
	\label{eq:forward_integral}
	f(\b x) 
	= 
	\int\limits_{\left[-\frac M2,\frac M2\right)^d} \hat f(\b v)\,\e^{2\pi\i \b v\b x}\,\mathrm d\b v,
	\quad \b x\in \R^d ,
\end{align}
from given measurements \mbox{$\hat f(\b k)\in\C$}, \mbox{$\b k\in\I_{\b{M}}$}, of its Fourier transform~$\hat f$.

To do so, this paper is organized as follows. 
Firstly, in Section~\ref{sec:nfft} we review the NFFT for trigonometric polynomials.
Subsequently, in Section~\ref{sec:shannon} we give an overview of the regularized Shannon sampling formulas, which play the key role in introducing the NFFT-like procedure for band\-limited functions in Section~\ref{sec:fast_algorithm_FT}.
Finally, in Section~\ref{sec:comparisonNFFT} we compare this new method to the classical NFFT.

\section{The NFFT\label{sec:nfft}}

For given nonequispaced nodes~\mbox{$\b x_j \in \T^d$},\, \mbox{$j=1,\dots,N$}, and given coefficients~\mbox{$\hat f_{\b k}\in\C$}, \mbox{$\b k\in\I_{\b{M}}$}, we consider the computation of the sums 
\begin{align}
	\label{eq:nfft}
	f(\b x_j) = \sum_{\b k \in \I_{\b{M}}} \hat{f}_{\b k}\, \e^{2\pi\i \b k \b x_j}, \quad j=1,\dots,N,
\end{align}
where the inner product shall be defined as usual as
\mbox{$\b k \b x \coloneqq k_1 x_1 + \dots + k_d x_d$}.
A fast approximate algorithm, the so-called \emph{nonequispaced fast Fourier transform~(NFFT)}, can be summarized as follows, see e.\,g.~\cite{duro93, bey95, st97, GrLe04, KeKuPo09} or~\cite[pp.~413--417]{PlPoStTa23}. \vspace{-2ex}

\begin{algorithm}{NFFT\vspace{0.5ex}}
	\label{alg:nfft}
	For~\mbox{$d,N \in \N$} and \mbox{$M \in 2\N$} let~\mbox{$\b x_j \in \T^d,\, j=1, \dots, N,$} be given nodes as well as~\mbox{$\hat f_{\b k} \in \C$}, \mbox{$\b k \in \I_{\b{M}}$}, given Fourier coefficients.
	Furthermore, we are given the oversampling factor~\mbox{$\sigma \geq 1$} with~\mbox{$M_{\sigma} \coloneqq 2 \lceil \,\lceil \sigma M \rceil / 2 \,\rceil \in 2\N$}, 
	as well as the window function~$\varphi$, the truncated function~$\varphi_m$ with~\mbox{$m \ll M_{\sigma}$}, and their \mbox{$1$-per}iodic versions~$\tilde \varphi$ and~\mbox{$\tilde \varphi_m$}.
	%
	\begin{enumerate}
		\item[0.] Precomputation:
		\begin{enumerate}
			\item Compute the nonzero Fourier coefficients~\mbox{$\hat \varphi(\b k)$} for \mbox{$\b k \in \I_{\b{M}}$}.
			\item Compute the values~\mbox{$\tilde\varphi_m \big(\b x_j-\frac{\b\ell}{M_{\b\sigma}}\big)$} for \mbox{$j = 1, \dots, N,$} as well as \mbox{$\b\ell\in \I_{\b{M_{\sigma}},m}(\b x_j)$}, cf.~\eqref{eq:indexset_x}.
		\end{enumerate}
		\item Set
		\hfill \mbox{$\mathcal O(|\I_{\b{M}}|)$}
		\begin{align*}
			\hat g_{\b k}\coloneqq
			\left\{
			\begin{array}{cl}
				\frac{\hat f_{\b k}}{\hat \varphi(\b k)} &\colon \b k \in \I_{\b{M}},\\
				0 &\colon \b k \in \I_{\b{M_{\b\sigma}}} \setminus \I_{\b{M}}.
			\end{array}
			\right.
		\end{align*}
		\item Compute
		\hfill \mbox{$\mathcal O(|\I_{\b{M}}|\log(|\I_{\b{M}}|))$}
		\begin{align*}
			g_{\b\ell}
			\coloneqq
			\frac{1}{|\I_{\b{M_{\b\sigma}}}|} \sum_{\b k \in \I_{\b{M}}}
			\hat g_{\b k}\, \e^{2\pi\i \b k \b\ell/M_{\sigma}},
			\quad \b\ell \in \I_{\b{M_{\b\sigma}}},
		\end{align*}
		by means of a \mbox{$d$-var}iate iFFT.
		\item Compute the short sums
		\hfill \mbox{$\mathcal O(N)$}
		\begin{align*}
			\tilde f_j 
			\coloneqq \hspace{-1em}
			\sum_{\b\ell \in \I_{\b{M_{\b\sigma}},m}(\b x_j)} \hspace{-1em} g_{\b\ell}\,\tilde \varphi_m \big(\b x_j-\tfrac{\b\ell}{M_{\sigma}}\big), \quad j = 1, \dots, N.
		\end{align*}
	\end{enumerate}
	\vspace{-1.1ex}
	\rule{\linewidth}{0.4pt}
	\textnormal{\textbf{Output:}} \mbox{$\tilde f_j \approx f(\b x_j)$} 
	\hfill
	\textnormal{\textbf{Complexity:}}~\mbox{$\mathcal O(|\I_{\b{M}}|\log(|\I_{\b{M}}|) + N)$} \hspace{-1.8ex} \vspace{0.5ex}
\end{algorithm}

\begin{Remark}
	\label{Rem:software}
	Note that Algorithm~\ref{alg:nfft} is part of the software packages~\cite{nfft3} and~\cite{FINUFFT}, respectively.
\end{Remark}

\noindent
By defining the \emph{nonequispaced Fourier matrix}
\begin{align*}
	\b A = \b A_{|\I_{\b{M}}|} \coloneqq \left( \e^{2\pi\i \b k \b x_j} \right)_{j=1,\,\b k \in \I_{\b{M}}}^{N} 
	\ \in \C ^{N\times |\I_{\b{M}}|},
\end{align*}
as well as the vectors
\mbox{$\b f\coloneqq\left(f(\b x_j)\right)_{j=1}^N$} and 
\mbox{$\b{\hat f}\coloneqq(\hat f_{\b k})_{\b k \in \I_{\b{M}}}$},
the computation of the sums in~\eqref{eq:nfft} can be written as
\mbox{$\b f = \b A \b{\hat f}$}.
By additionally defining the diagonal matrix
\begin{align}
	\label{eq:matrix_D}
	\b D 
	\coloneqq 
	\text{diag} \left( \frac 1{|\I_{\b{M_{\b\sigma}}}|\cdot\hat{\varphi}(\b k)} \right)_{\b k \in \I_{\b{M}}} 
	\in \C^{|\I_{\b{M}}|\times |\I_{\b{M}}|},
\end{align}
the truncated \emph{Fourier matrix}
\begin{align}
	\label{eq:matrix_F}
	\b F 
	\coloneqq 
	\left( \e^{2\pi\i \b k \b\ell/M_{\b\sigma}} \right)_{\b\ell \in \I_{\b{M_{\b\sigma}}},\, \b k \in \I_{\b{M}}} 
	\in \C ^{|\I_{\b{M_{\sigma}}}|\times |\I_{\b{M}}|},
\end{align}
and the \mbox{$(2m+1)^d$}-sparse matrix
\begin{align}
	\label{eq:matrix_B}
	\b B 
	\coloneqq 
	\bigg( \tilde \varphi_m \big(\b x_j-\tfrac{\b\ell}{M_{\sigma}}\big) \bigg)_{j=1,\, \b\ell \in \I_{\b{M_{\b\sigma}}}}^{N} 
	\in \R^{N\times |\I_{\b{M_{\b\sigma}}}|},
\end{align}
where by definition of the index set 
\begin{align}
	\label{eq:indexset_x}
	\I_{\b{M_{\b\sigma}},m}(\b x_j)
	&\coloneqq
	\left\{ \b\ell \in \I_{\b{M_{\b\sigma}}} \colon \exists \b z \in \Z^d \text{ with } \right. 
	\left. -m \cdot \b 1_d \leq M_{\sigma} \cdot (\b x_j + \b z) - \b\ell \leq m \cdot \b 1_d \right\} 
\end{align}
each row of~$\b B$ contains at most \mbox{$(2m+1)^d$} nonzeros, the NFFT in Algorithm~\ref{alg:nfft} can be formulated in matrix-vector notation as
\mbox{$\b A \approx \b B \b F \b D$}, cf.~\cite[p.~419]{PlPoStTa23}.
This is to say, using the definition of these matrices, the NFFT performs the approximation
\begin{align}
	\label{eq:approx_nfft}
	\e^{2\pi\i \b k \b x_j} 
	\approx 
	\hspace{-0.5em} \sum_{\b\ell \in \I_{\b{M_{\b\sigma}},m}(\b x_j)} \hspace{-1em}
	\frac{\e^{2\pi\i \b k \b\ell/M_{\sigma}} \,\tilde \varphi_m \big(\b x_j-\tfrac{\b\ell}{M_{\sigma}}\big)}{|\I_{\b{M_{\b\sigma}}}|\cdot\hat{\varphi}(\b k)}
\end{align}
for~\mbox{$\b k\in\I_{\b M}$} and~\mbox{$\b x_j\in\T^d$}, \mbox{$j=1,\dots,N$}.


\section{Regularized Shannon sampling formulas \label{sec:shannon}}

A function~\mbox{$f \colon \R^d\to\C$} is said to be \emph{band\-limited} with \emph{band\-width}~\mbox{$M\in\N$}, if
the support of its \emph{(continuous) Fourier transform}
\begin{align}
	\label{eq:inverse_integral}
	\hat f(\b v)
	\coloneqq 
	\int\limits_{\R^d} 
	f(\b x)\,\e^{-2\pi\i \b v\b x}\,\mathrm d\b x,
	\quad \b v\in\R^d,
\end{align}
is contained in~\mbox{$\left[-\frac M2,\frac M2\right]^d$}.
The space of all band\-limited functions with band\-width~\mbox{$M$} shall be denoted by
\begin{align*}
	{\mathcal B}_{M/2}(\R^d) \coloneqq \Big\{ f \in L_2(\R^d)
	\colon\, \mathrm{supp}(\hat f)
	\!\subseteq\! \left[- \tfrac{M}{2},\,\tfrac{M}{2}\right]^d \!\Big\} , \!
\end{align*}
which is also known as the \emph{Paley--Wiener space}.
Note that
\begin{align}
	\label{eq:embedding_BM2}
	{\mathcal B}_{M/2}(\R^d) \subseteq L_2(\R^d) \cap C_0(\R^d) \cap C^\infty(\R^d) ,
\end{align}
cf.~\cite[Lemma~4.1]{Kircheisdiss}.
Thus, the Fourier inversion theorem, see e.\,g.~\cite[Theorem~2.23]{PlPoStTa23}, guarantees that the \emph{inverse Fourier transform} of~$f$ can be written as given in~\eqref{eq:forward_integral}.

By the famous Whittaker--Kotelnikov--Shannon sampling theorem (\cite{Whittaker, Kotelnikov, Shannon49}) any band\-limited function~\mbox{$f\in {\mathcal B}_{M/2}(\R^d)$} can be recovered from its samples~\mbox{$f\big(\tfrac{\b\ell}{L}\big)$}, \mbox{$\b\ell\in\Z^d$}, with~\mbox{$L \geq M$}, \mbox{$L\in\N$}, in the form
\begin{align}
	\label{eq:sampling_theorem}
	f(\b x) 
	&=
	\sum_{\b\ell\in\Z^d} 
	f\big(\tfrac{\b\ell}{L}\big) \,	\sinc \big(L\pi\big(\b x - \tfrac{\b\ell}{L}\big)\big),
	\quad \b x\in\R^d , 
\end{align}
where the \mbox{$\sinc$ func}tion is given by~\mbox{$\sinc(\b x) \coloneqq \prod_{t=1}^d \sinc(x_t)$} with
\begin{align*}
	\sinc(x) \coloneqq \left\{ \begin{array}{ll}  \frac{\sin x}{x} & \colon x \in \R \setminus \{0\} , \\ [1ex]
		1 & \colon x = 0 . \end{array} \right.
\end{align*}
It is well known that the series in~\eqref{eq:sampling_theorem} converges absolutely and uniformly on whole~\mbox{$\R^d$}.

Unfortunately, the numerical use of this classical Whittaker--Kotelnikov--Shannon sampling series~\eqref{eq:sampling_theorem} is limited, since it requires infinitely many samples, which is impossible in practice, and its truncated version is not a good approximation due to the slow decay of the \mbox{$\sinc$ func}tion, see~\cite{Ja66}.
In addition to this rather poor convergence, it is known, see~\cite{Fe92a,Fe92b,DDeV03}, that in the presence of noise in the samples~\mbox{$f\big(\tfrac{\b\ell}{L}\big)$}, \mbox{$\b\ell\in\Z^d$}, the convergence of the Shannon sampling series~\eqref{eq:sampling_theorem} may even break down completely.

Based on this observation, numerous approaches for numerical realizations have been developed, where the Shannon sampling series was regularized with a suitable window function.
Note that many authors such as~\cite{dau92, Nat86, Rap96, Par97, StTa06} used window functions in the frequency domain, but the recent study~\cite{KiPoTa23} has shown that it is much more beneficial to employ a window function in the spatial domain, cf.~\cite{Q03, Q04, StTa06, MXZ09, LZ16, ChZh19, KiPoTa22}.

Therefore, for a given \mbox{$m\in\N$} with \mbox{$2m \ll L$} we introduce the set $\Phi_{m,L}$ of all window functions \mbox{$\varphi \colon\, \R \to [0,\,1]$} with the following properties: 
\begin{itemize}
	\item $\varphi$ is compactly supported on~\mbox{$\left[-\tfrac{m}{L},\,\tfrac{m}{L}\right]$}, belongs to $L_1(\R) \cap C_0(\R)$ and is even.
	\item $\varphi$ restricted to ${[0,\,\infty)}$ is monotonously non-increasing with \mbox{$\varphi(0) = 1$}. 
\end{itemize}
%

\begin{Remark}
	\label{Rem:window_functions_Shannon}
	As examples of such window functions we consider 
	the {\emph{\mbox{$\sinh$-type} window function}}
	\begin{align}
		\label{eq:varphisinh}
		\!\!\varphi_{\sinh}(x) \coloneqq
			\frac{1}{\sinh \beta}\, \sinh\!\Big(\beta\sqrt{1-\big(\tfrac{Lx}{m}\big)^2}\,\Big)
			\, \chi_{\left[-\tfrac{m}{L},\,\tfrac{m}{L}\right]}(x) 
		\!
	\end{align}
	with certain~\mbox{$\beta > 0$},
	and the \emph{continuous Kaiser--Bessel window function}
	\begin{align*}
		\varphi_{\mathrm{cKB}}(x) \coloneqq
			\frac{\bigg( I_0\Big(\beta\,\sqrt{1-\big(\tfrac{Lx}{m}\big)^2}\,\Big) - 1\bigg)}{I_0(\beta) - 1}\, 
			\, \chi_{\left[-\tfrac{m}{L},\,\tfrac{m}{L}\right]}(x) 
	\end{align*}
	with certain~\mbox{$\beta > 0$}, where~\mbox{$I_0$} denotes the \emph{modified Bessel function of first kind}. 
	Note that these window functions are well-studied in the context of the NFFT, see e.\,g.~\cite{PoTa21a}.
\end{Remark}

Then, for a fixed window function~\mbox{$\varphi \in \Phi_{m,L}$} we study the \emph{regularized Shannon sampling formula with localized sampling}
\begin{align}
	\label{eq:Rmf(x)_multi}
	(R_{\varphi,m}f)({\b x}) 
	&\coloneqq
	\sum_{{\b \ell} \in \Z^d} f\big(\tfrac{{\b \ell}}{L}\big) \,\sinc\big(L\pi \big(\b x - \tfrac{\b \ell}{L}\big)\big) \, \varphi\big({\b x} - \tfrac{{\b \ell}}{L}\big) , \quad \b x \in \R^d .
\end{align}
Note that this is an \emph{interpolating approximation} of~$f$ on \mbox{$\frac{1}{L}\,\Z^d$}, i.\,e., we have 
\begin{align*}
	f\big(\tfrac{\b k}{L}\big) = (R_{\varphi,m}f)\big(\tfrac{\b k}{L}\big), \quad \b k \in \Z^d ,
\end{align*}
since by assumption~\mbox{$\varphi(0) = 1$} and~\mbox{$\sinc(\pi(k-\ell)) = \delta_{k,\ell}$} for all~\mbox{$k,\,\ell\in \Z$} with the Kronecker symbol~\mbox{$\delta_{k,\ell}$}.
Then it is known that the regularized Shannon sampling formula~\mbox{$R_{\varphi,m}f$} in~\eqref{eq:Rmf(x)_multi} with suitable window function~\mbox{$\varphi\in\Phi_{m,L}$} yields a good approximation of~$f$, cf.~\cite{KiPoTa22,KiPoTa23,Kircheisdiss}.

\section{NFFT-like procedure for bandlimited functions \label{sec:fast_algorithm_FT}}

Now assume we are given the values~\mbox{$\hat{f}(\b k)$}, \mbox{$\b k\in \I_{\b M}$}, of the Fourier transform~\eqref{eq:forward_integral} of a band\-limited function~\mbox{$f\in\mathcal{B}_{M/2}(\R^d)$}, and
we are looking for function evaluations~\mbox{$f(\b x_j)$} at given nonequispaced points~$\b x_j$, \mbox{$j=1,\dots N$}.
Further we assume that the function~\mbox{$f$} fulfills the condition
\begin{equation}
	\label{eq:sampledecay}
	\sum_{\b\ell\in \mathbb Z^d} \big| f\big(\tfrac{\b\ell}{L}\big) \big| < \infty ,
\end{equation}
such that the Fourier series
\begin{equation*}
	\hat f(\b v) = \frac{1}{L^d} \sum_{\b\ell\in \mathbb Z^d} f\big(\tfrac{\b\ell}{L}\big) \,\e^{-2\pi\i \b\ell \b v/L} ,
	\quad \b v\in\left[-\tfrac{L}{2},\,\tfrac{L}{2}\right]^d ,
\end{equation*}
converges absolutely and uniformly, see~\cite[p.~12]{KiPoTa23}, and thus the considered problem is well-defined.
Note that by~\cite[Lemma~2]{ScSi00} it is known that~\mbox{$f\in L_1(\R^d)$} fulfill the condition~\eqref{eq:sampledecay}.
Moreover, \mbox{$f\in L_1(\R^d)$} directly implies~\mbox{$\hat f\in C_0(\R^d)$}.

\begin{Remark}
	Note that the recent work~\cite{EhGrKl24} derived error estimates for a familiar problem, where however for functions~\mbox{$f\in C(\R)$} satisfying certain decay and smoothness conditions and equispaced points~\mbox{$\b x_j=\frac{2j- N}{2N}$}, \mbox{$j=1,\dots,N$}, simply the FFT can be used.
\end{Remark}

In order to compute the values~\mbox{$f(\b x_j)$}, \mbox{$j=1,\dots N$}, we aim to make use of the regularized Shannon sampling formulas, see Section~\ref{sec:shannon}.
Inserting the approximation~\eqref{eq:Rmf(x)_multi} into the Fourier transform~\eqref{eq:inverse_integral} and using the definition of the regularized \mbox{$\sinc$ func}tion
\begin{align}
\label{eq:sinc_reg}
	\psi(\b x) \coloneqq \sinc(L\pi \b x) \,\varphi(\b x),
\end{align}
we have
\begin{align}
	\label{eq:approx_fourier_trafo}
	\hat f(\b v)
	&=	
	\int_{\R^d} f(\b x)\,\e^{-2\pi\i \b v \b x}\,\mathrm d\b x
	\approx
	\int_{\R^d} (R_{\varphi,m} f)(\b x)\,\e^{-2\pi\i \b v \b x}\,\mathrm d\b x \notag \\
	&=	
	\int_{\R^d} \,\sum_{\b\ell\in \Z^d} f\big(\tfrac{\b\ell}{L}\big) \, \psi\big(\b x-\tfrac{\b\ell}{L}\big)\,\e^{-2\pi\i \b v \b x}\,\mathrm d\b x \notag \\
	&=
	\sum_{\b\ell\in \Z^d} f\big(\tfrac{\b\ell}{L}\big) \,\e^{-2\pi\i \b v \b\ell/L} \int_{\R^d} \psi(\b y)\,\e^{-2\pi\i \b v \b y}\,\mathrm d\b y \notag \\ &
	=
	\bigg(\sum_{\b\ell\in \Z^d} f\big(\tfrac{\b\ell}{L}\big) \,\e^{-2\pi\i \b v \b\ell/L}\bigg) \cdot \hat\psi(\b v) ,
\end{align}
where summation and integration may be interchanged 
by the theorem of Fubini--Tonelli.
%
By defining
\begin{align}
	\label{eq:def_nuhat}
	\hat\nu(\b v) \coloneqq \sum_{\b\ell\in \Z^d} f\big(\tfrac{\b\ell}{L}\big) \,\e^{-2\pi\i \b v \b\ell/L} ,
	\quad \b v\in\R^d ,
\end{align}
we recognize that this function~\mbox{$\hat\nu$} is \mbox{$L$-per}iodic. 
Thus, due to the fact that the Fourier transform of the band\-limited function~\mbox{$f\in\mathcal{B}_{M/2}(\R^d)$} is non-periodic, the approximation~\eqref{eq:approx_fourier_trafo} can only be reasonable for~\mbox{$\b v\in\big[-\frac L2,\frac L2\big]^d$}.

As the goal is to recover the nonequispaced samples~\mbox{$f(\b x_j)$}, \mbox{$j=1,\dots,N$}, by means of a regularized Shannon sampling formula~\eqref{eq:Rmf(x)_multi}, we need access to as many equispaced samples~\mbox{$f\big(\tfrac{\b\ell}{L}\big)$} as possible, i.\,e., we are looking for an inversion formula for~\eqref{eq:def_nuhat}.
To this end, note that~\eqref{eq:def_nuhat} can be written as
\begin{align*}
	\hat\nu(\b v)
	&= 
	\sum_{\b\ell\in\I_{\b\Theta}} f\big(\tfrac{\b\ell}{L}\big) \,\e^{-2\pi\i \b v \b\ell/L} 
	+ 
	\sum_{\b r \in \Z^d\setminus\{\b 0\}} \sum_{\b\ell\in\I_{\b\Theta}} f\big(\tfrac{\b\ell+\b r\Theta}{L}\big)\, \e^{-2\pi\i \b v (\b\ell+\b r\Theta)/L} ,
	\quad \b v\in\R^d ,
\end{align*}
with the index set~\mbox{$\I_{\b\Theta}$} with~\mbox{$\b \Theta=\Theta\cdot \b 1_d$}, \mbox{$\Theta\in 2\N$}.
Since~\mbox{$f\in {\mathcal B}_{M/2}(\R^d) \subseteq C_0(\R^d)$}, see~\eqref{eq:embedding_BM2}, the equispaced samples~\mbox{$f\big(\tfrac{\b\ell}{L}\big)$} are negligible for all~\mbox{$\|\b\ell\|_\infty \geq \frac{\Theta}{2}$} with suitably chosen~$\Theta$.
In order to avoid aliasing in the computation we assume that~\mbox{$\Theta=L$} is sufficient.
Hence, we consider
\begin{align}
	\label{eq:def_thetahat}
	\hat\nu(\b v)
	\approx
	\hat\vartheta(\b v)
	&\coloneqq 
	\sum_{\b\ell\in\I_{\b L}} f\big(\tfrac{\b\ell}{L}\big) \,\e^{-2\pi\i \b v \b\ell/L} ,
	\quad \b v\in\R^d ,
\end{align}
and thus by~\eqref{eq:approx_fourier_trafo} the approximation
\begin{align}
	\label{eq:approx_fourier_trafo_theta}
	\hat f(\b v)
	&\approx
	\hat\vartheta(\b v) \cdot \hat\psi(\b v) ,
	\quad \b v\in\big[-\tfrac L2,\tfrac L2\big]^d . 
\end{align}
Since it is additionally known that \mbox{$\hat f(\b v)=0$} for all \mbox{$\b v\notin\big[-\frac M2,\frac M2\big]^d$} and \mbox{$\hat\psi(\b v)\neq 0$} for all~\mbox{$\b v\in\big[-\frac L2,\frac L2\big]^d$}, we might use~\eqref{eq:approx_fourier_trafo_theta} and~\eqref{eq:def_thetahat} for given \mbox{$\hat{f}(\b k)$}, \mbox{$\b k\in \I_{\b M}$}, to approximate the equispaced samples~\mbox{$f\big(\tfrac{\b\ell}{L}\big)$}, \mbox{$\b\ell\in\I_{\b L}$}, by setting
\begin{align*}
	\hat\vartheta(\b k)
	=
	\left\{
	\begin{array}{cl}
		\frac{\hat f(\b k)}{\hat\psi(\b k)} &\colon\quad \b k\in\I_{\b M} , \\
		0&\colon\quad \b k\in\I_{\b L}\setminus\I_{\b M} ,
	\end{array}
	\right.
\end{align*}
and subsequently computing
\begin{align}
	\label{eq:def_theta}
	f\big(\tfrac{\b\ell}{L}\big)
	\approx
	\vartheta_{\b\ell}
	\coloneqq
	\frac{1}{|\I_{\b L}|} \sum_{\b k\in\I_{\b L}} \hat\vartheta(\b k) \,\e^{2\pi\i \b k \b \ell/L} ,
	\quad \b\ell\in\I_{\b L} ,
\end{align}
by means of an iFFT.

To finally approximate the samples~\mbox{$f(\b x_j)$}, \mbox{$j=1,\dots,N$}, we make use of the regu\-la\-rized Shannon sampling formula~\eqref{eq:Rmf(x)_multi}.
Note that since we assumed that the window function~\mbox{$\varphi\in\Phi_{m,L}$} is compactly supported, the computation of~\mbox{$(R_{\varphi,m}f)(\b x)$} for fixed~\mbox{$\b x \in \R^d \setminus \tfrac{1}{L}\,\Z^d$} requires only~\mbox{$(2m+1)^d$} samples~\mbox{$f\big(\tfrac{\b \ell}{L}\big)$}. 
However, we have already encountered that~\eqref{eq:def_theta} can only be used to approximate~\mbox{$f\big(\tfrac{\b\ell}{L}\big)$} for~\mbox{$\b\ell\in\I_{\b L}$} in order to avoid aliasing in the computation of the inverse Fourier transform in~\eqref{eq:def_theta}.
Thereby, we are confronted with a limitation of the feasible points to~\mbox{$\b x_j \in [-\tfrac 12+\tfrac mL, \tfrac 12-\tfrac mL)^d$}, \mbox{$j=1,\dots,N$}, since only in this case exclusively the evaluations \mbox{$f\big(\tfrac{\b\ell}{L}\big)$}, \mbox{$\b\ell\in\I_{\b L}$}, are needed for the computation. 
Hence, the final approximation is computed by
\begin{align*}
	(R_{\varphi,m} f)(\b x_j)
	\approx
	f_j
	&\coloneqq
	\sum_{\b\ell \in \I_{\b L}} \vartheta_{\b\ell} \,\psi \big(\b x_j - \tfrac{\b\ell}{L}\big) 
	=
	\sum_{\b\ell \in \mathcal J_{\b{L},m}(\b x_j)} \vartheta_{\b\ell} \,\psi \big(\b x_j - \tfrac{\b\ell}{L}\big) , \notag
\end{align*}
where the index set of the nonzero entries
\begin{align}
	\label{eq:indexset_x_nonperiodic}
	\mathcal J_{\b{L},m}(\b x_j)
	\coloneqq
	\left\{ \b\ell\in\Z^d \colon -m+L \b x_j \leq \b\ell \leq m+L \b x_j \right\}
\end{align}
contains at most \mbox{$(2m+1)^d$} entries for each fixed~$\b x_j$, cf.~\eqref{eq:indexset_x}.
Thus, the obtained algorithm can be summarized as follows, cf.~\cite[Algorithm~5.16]{Kircheisdiss}.

\begin{algorithm}{NFFT-like procedure for bandlimited functions}
	\label{alg:nfft_generalized}
	For~\mbox{$d,m,N \in \N$}, \mbox{$M\in 2\N$}, and~\mbox{$L=M(1+\lambda)\in\N$} with oversampling parameter~\mbox{$\lambda\geq 0$} let~\mbox{$\b x_j \in [-\tfrac 12+\tfrac mL, \tfrac 12-\tfrac mL)^d$}, {$j=1,\dots,N$}, be given nodes as well as~\mbox{$\hat f(\b k) \in \C$}, \mbox{$\b k \in \I_{\b{M}}$}, given evaluations of the Fourier transform of the band\-limited function~\mbox{$f\in\mathcal{B}_{M/2}(\R^d)$}.
	Furthermore, we are given the window function~\mbox{$\varphi\in\Phi_{m,L}$}, the corresponding regularized \mbox{$\sinc$ func}tion~$\psi$ in~\eqref{eq:sinc_reg}, and its Fourier transform~\mbox{$\hat\psi$}.
	\begin{enumerate}
		\item[0.] Precomputation:
		\begin{enumerate}
			\item Compute the nonzero values \mbox{$\hat\psi(\b k)$} for \mbox{$\b k \in \I_{\b{M}}$}.
			\item Compute the evaluations~\mbox{$\psi\big(\b x_j-\tfrac{\b\ell}{L}\big)$} for \mbox{$j = 1, \dots, N,$} as well as \mbox{$\b\ell\in \mathcal J_{\b{L},m}(\b x_j)$}, cf.~\eqref{eq:indexset_x_nonperiodic}.
		\end{enumerate}
		\item Set \hfill \mbox{$\mathcal O(|\I_{\b M}|)$}
		\begin{align*}
			\hat\vartheta(\b k) \coloneqq 
			\left\{
			\begin{array}{cl}
				\frac{\hat f(\b k)}{\hat\psi(\b k)} &\colon\quad \b k\in\I_{\b M} , \\
				0&\colon\quad \b k\in\I_{\b L}\setminus\I_{\b M}.
			\end{array}
			\right.
		\end{align*}
		\item Compute
		\hfill \mbox{$\mathcal O(|\I_{\b{M}}|\log(|\I_{\b{M}}|))$}
		\begin{align*}
			\vartheta_{\b\ell} \coloneqq \frac{1}{|\I_{\b L}|} \sum_{\b k\in \I_{\b L}} \hat\vartheta(\b k) \,\e^{2\pi\i \b k \b\ell/L}, \quad \b\ell\in\I_{\b L},
		\end{align*}
		by means of a \mbox{$d$-var}iate iFFT.
		\item Compute the short sums \hfill \mbox{$\mathcal O(N)$}
		\begin{align*}
			f_j \coloneqq \sum_{\b\ell \in \mathcal J_{\b{L},m}(\b x_j)} \vartheta_{\b\ell}\,\psi\big(\b x_j-\tfrac{\b\ell}{L}\big), \quad j=1,\dots, N.
		\end{align*}
	\end{enumerate}
	\vspace{-1.1ex}
	\rule{\linewidth}{0.4pt}
	\textnormal{\textbf{Output:}} \mbox{$f_j \approx f(\b x_j)$} 
	\hfill
	\textnormal{\textbf{Complexity:}} \mbox{$\mathcal O(|\I_{\b{M}}|\log(|\I_{\b{M}}|) + N)$} \hspace{-1.8ex} \vspace{0.5ex}
\end{algorithm}

Note that by defining the vector~\mbox{$\b{\hat f} \coloneqq ( \hat f(\b k) )_{\b k\in\I_{\b M}}$} as well as the diagonal matrix
\begin{align}
	\label{eq:matrix_Dpsi}
	\b D_{\hat\psi}
	\coloneqq 
	\text{diag} \left( \frac 1{|\I_{\b{L}}|\cdot\hat{\psi}(\b k)} \right)_{\b k \in \I_{\b{M}}} 
	\ \in \C^{|\I_{\b{M}}|\times |\I_{\b{M}}|} 
\end{align}
and the \mbox{${(2m+1)^d}$-sparse} matrix
\begin{align}
	\label{eq:matrix_Psi}
	\b\Psi
	\coloneqq
	\bigg( \psi\big(\b x_j - \tfrac{\b\ell}{L}\big) \bigg)_{j=1,\, \b\ell\in\I_{\b L}}^{N}
	\ \in \R^{N\times |\I_{\b L}|} ,
\end{align}
the approximation of Algorithm~\ref{alg:nfft_generalized} is given by
\begin{align}
	\label{eq:function_eval}
	\b f 
	= \b\Psi \b F \b D_{\hat\psi} \b{\hat f} ,
\end{align}
where~\mbox{$\b F\in \C^{|\I_{\b L}|\times |\I_{\b M}|}$} denotes the Fourier matrix~\eqref{eq:matrix_F} with~\mbox{$L=M_\sigma$}.

\section{Comparison to the classical NFFT \label{sec:comparisonNFFT}}

Note that one might also directly apply an equispaced quadrature rule to the inverse Fourier transform~\eqref{eq:forward_integral}, i.\,e., consider the approximation
\begin{align*}
	f(\b x) 
	= \int\limits_{[-\frac M2,\frac M2]^d} \hat f(\b v) \,\e^{2\pi\i \b v \b x} \,\mathrm d\b v
	\approx \sum_{\b k\in \I_{\b M}} \hat f(\b k) \,\e^{2\pi\i \b k \b x} ,
\end{align*}
such that the function evaluations~\mbox{$f(\b x_j)$}, \mbox{$j=1,\dots,N$}, could also be approximated efficiently by means of an NFFT. 
Since this raises the question of which of the two methods, Algorithm~\ref{alg:nfft} or Algorithm~\ref{alg:nfft_generalized}, is more advantageous, this section deals with the comparison of the two approaches.

Considering the matrix notations~\mbox{$\b B \b F \b D$} and~\mbox{$\b\Psi \b F \b D_{\hat{\psi}}$}, cf.~\eqref{eq:approx_nfft} and~\eqref{eq:function_eval}, the first thing to realize is that for~\mbox{$\b B\in \R^{N\times |\I_{\b L}|}$} in~\eqref{eq:matrix_B} the window function~\mbox{$\varphi_m(\b x)$} is used, while for~\mbox{$\b\Psi\in \R^{N\times |\I_{\b L}|}$} in~\eqref{eq:matrix_Psi} we consider the regularized \mbox{$\sinc$ func}tion~\mbox{$\psi(\b x)$} in~\eqref{eq:sinc_reg}.
A similar remark can also be made about the diagonal matrices~\mbox{$\b D\in \C^{|\I_{\b{M}}|\times |\I_{\b{M}}|}$} in~\eqref{eq:matrix_D} and~\mbox{$\b D_{\hat{\psi}}\in \C^{|\I_{\b{M}}|\times |\I_{\b{M}}|}$} in~\eqref{eq:matrix_Dpsi}.

Additionally, it is important to note that the two methods can only be compared for~\mbox{$\b x\in[-\tfrac 12+\tfrac mL, \tfrac 12-\tfrac mL)^d$}, as the approximation by Algorithm~\ref{alg:nfft_generalized} is only reasonable in this case.
This implies that the matrix~$\b B$ in~\eqref{eq:matrix_B} is, unlike usual, non-periodic, whereas the matrix~$\b\Psi$ in~\eqref{eq:matrix_Psi} is inherently non-periodic by definition.

To study the quality of both approaches, note that by the NFFT we are given the approximation
\begin{align}
	\label{eq:approx_exp_nfft}
	\e^{2\pi\i \b k \b x} 
	&\approx 
	\frac 1{|\I_{\b{L}}|\cdot\hat{\varphi}(\b k)} 
	\sum_{\b\ell \in \I_{\b{L}}} 
	\e^{2\pi\i \b k \b\ell/L}
	\,\tilde \varphi_m \hspace{-2.5pt}\left(\b x-\tfrac{\b\ell}{L}\right) , \quad \b x \in \T^d ,
\end{align}
for~\mbox{$\b k\in\I_{\b M}$} fixed, cf.~\eqref{eq:approx_nfft} with~\mbox{$L=M_\sigma$}, where~\mbox{$\tilde \varphi_m (\b x) = \sum_{\b r\in\Z^d} \varphi_m(\b x+\b r)$} denotes the \mbox{$1$-per}iodic version of the compactly supported window function~\mbox{$\varphi_m$}.
Thus, we look for a comparable approximation of the exponential function using our newly proposed method in Algorithm~\ref{alg:nfft_generalized}.
For this purpose, note that~\mbox{$g(\b x) \coloneqq \hat\psi(\b x)\,\e^{2\pi\i\b k \b x}$} with~\mbox{$\b k\in\R^d$} fixed possesses the Fourier transform~\mbox{$\hat g(\b v) = \psi(\b k - \b v)$}.
Therefore, we have~\mbox{$g\in\mathcal{B}_{M/2}(\R^d)$} for all~\mbox{$\b k\in\big[-\frac M2+\frac mL,\frac M2-\frac mL\big]^d$}, i.\,e., considering~\eqref{eq:Rmf(x)_multi} for this function~$g$ yields 
\begin{align*}
	\hat\psi(\b x)\,\e^{2\pi\i\b k \b x} 
	\approx
	\sum_{{\b \ell} \in \Z^d} \hat\psi\big(\tfrac{{\b \ell}}{L}\big)\,\e^{2\pi\i\b k \b\ell/L} \,\psi\big({\b x} - \tfrac{{\b \ell}}{L}\big) ,
	\quad \b x \in \R^d ,
\end{align*}
or rather
\begin{align*}
	\e^{2\pi\i\b k \b x} 
	\approx
	\sum_{{\b \ell} \in \I_{\b L}} \frac{\hat\psi\big(\tfrac{{\b \ell}}{L}\big)}{\hat\psi(\b x)}\,\e^{2\pi\i\b k \b\ell/L} \,\psi\big({\b x} - \tfrac{{\b \ell}}{L}\big) ,
	\quad \b x \in \big[-\tfrac 12+\tfrac mL, \tfrac 12-\tfrac mL\big)^d .
\end{align*}
Since numerical experiments have shown that~\mbox{$\hat{\psi}(\b y) \approx \frac{1}{|\I_{\b L}|}$}, \mbox{$\b y \in \big[-\tfrac M2, \tfrac M2\big)^d$}, for the window functions mentioned in Remark~\ref{Rem:window_functions_Shannon}, the above approximation simplifies to
\begin{align}
	\label{eq:approx_exp_nfft_like}
	\e^{2\pi\i\b k \b x} 
	\approx
	\sum_{{\b \ell} \in \I_{\b L}} \e^{2\pi\i\b k \b\ell/L} \,\psi\big({\b x} - \tfrac{{\b \ell}}{L}\big) ,
	\quad \b x \in \big[-\tfrac 12+\tfrac mL, \tfrac 12-\tfrac mL\big)^d ,
\end{align}
which equals the approximation~\mbox{$\b\Psi \b F \b D_{\hat\psi}$} of Algorithm~\ref{alg:nfft_generalized}, since~\mbox{$|\I_{\b L}|\,\hat{\psi}(\b k) \approx 1$}, \mbox{$\b k\in\I_{\b M}$}.
Therefore, we can compare the quality of the two methods by considering the approximations~\eqref{eq:approx_exp_nfft} and~\eqref{eq:approx_exp_nfft_like} of the exponential function.

For simplicity we restrict ourselves to the one-dimensional setting~\mbox{$d=1$} for the visualization.
To estimate the quality of the approaches, we consider the approximation error
\begin{align}
	\label{eq:maxerr_comparison_nfft}
	e(v)
	\coloneqq
	\max_{x_p, p=1,\dots,P} \big|\mathrm e^{2\pi\mathrm i v x_p} - h(x_p)\big| ,
\end{align}
where the term~\mbox{$h(x_p)$} is a placeholder for the right-hand sides of~\eqref{eq:approx_exp_nfft} and~\eqref{eq:approx_exp_nfft_like}, respectively, evaluated at a fine grid of~\mbox{$P=10^5$} equispaced points~\mbox{$x_p$}, \mbox{$p=1,\dots,P$}.
This approximation error~\eqref{eq:maxerr_comparison_nfft} shall now be computed for several values 
\begin{align}
	\label{eq:maxerr_comparison_nfft_evaluation_points}
	\hspace{2em}
	v_s &= -\tfrac M2-m+\tfrac{s}{S} \in \big[-\tfrac M2-m,\tfrac M2+m\big] , 
	\quad s=0,\dots,S(M+2m) , 
\end{align}
where~\mbox{$S=1$} corresponds to integer evaluation, whereas we use~\mbox{$S=32$} to examine the approximation at non-integer points as well.
Note that~\eqref{eq:approx_exp_nfft} is expected to provide a good approximation only for~\mbox{$v\in\big[-\frac M2,\frac M2\big]$}, while~\eqref{eq:approx_exp_nfft_like} is expected to do so only for~\mbox{$v\in\big[-\frac M2+\frac mL,\frac M2-\frac mL\big]$}.
Nevertheless, we test for~$v$ from a larger interval to confirm these assumptions.

\begin{figure}[ht]
	\centering
	\captionsetup[subfigure]{justification=centering}
	\begin{subfigure}[t]{0.4\columnwidth}
		\includegraphics[width=\textwidth]{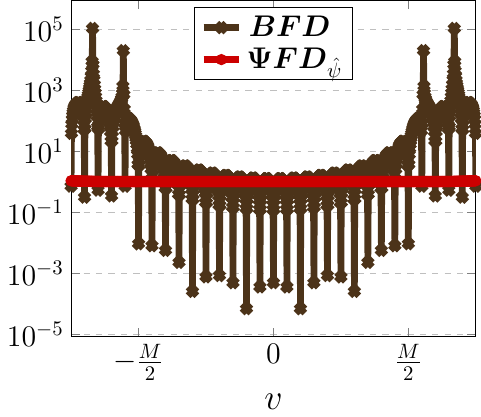}
		\caption{\mbox{$x \in \big[-\tfrac 12, \tfrac 12\big)$}}
		\label{fig:comparison_nfft_full_interval}
	\end{subfigure}
	\hspace{1em}
	\begin{subfigure}[t]{0.4\columnwidth}
		\includegraphics[width=\textwidth]{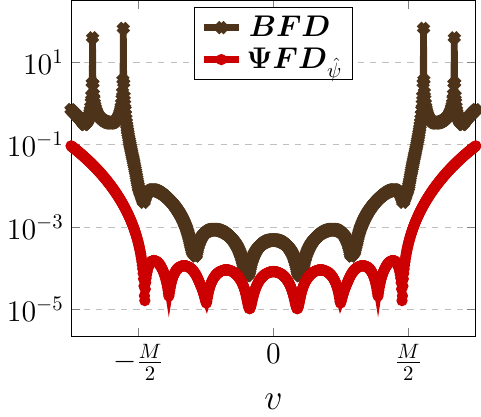}
		\caption{\mbox{$x \in \big[-\tfrac 12+\tfrac mL, \tfrac 12-\tfrac mL\big)$\hspace{-1em}}}
		\label{fig:comparison_nfft_truncated_interval}
	\end{subfigure}
	\caption{Maximum approximation error~\eqref{eq:maxerr_comparison_nfft} for~\mbox{$P=10^5$} computed for~\eqref{eq:maxerr_comparison_nfft_evaluation_points} with~\mbox{$S=32$} using the $\sinh$-type window function~\eqref{eq:varphisinh} as well as~\mbox{$M=20$}, \mbox{$\lambda=1$}, \mbox{$L=(1+\lambda)M$}, and~\mbox{$m=5$} in the one-dimensional setting~\mbox{$d=1$}.
		\label{fig:comparison_nfft}}
\end{figure}

The corresponding outcomes when computing the approximations~\eqref{eq:approx_exp_nfft} and~\eqref{eq:approx_exp_nfft_like} using the $\sinh$-type window function~\eqref{eq:varphisinh} as well as the parameters~\mbox{$M=20$}, \mbox{$\lambda=1$}, \mbox{$L=(1+\lambda)M$}, and~\mbox{$m=5$}, are displayed in Figure~\ref{fig:comparison_nfft}.
For~\mbox{$x \in \big[-\tfrac 12, \tfrac 12\big)$} it is easy to see that our newly proposed method~\eqref{eq:approx_exp_nfft_like} indeed does not provide reasonable results, while the approximation~\eqref{eq:approx_exp_nfft} by means of the NFFT is only useful at integer points~$v$.
For the truncated interval~\mbox{$x \in \big[-\tfrac 12+\tfrac mL, \tfrac 12-\tfrac mL\big)$}, however, both approximations~\eqref{eq:approx_exp_nfft} and~\eqref{eq:approx_exp_nfft_like} are clearly beneficial for non-integer points~$v$ as well, but as expected these methods only succeed when~\mbox{$|v|\leq \frac M2$}.
Nevertheless, although also the approximation~\eqref{eq:approx_exp_nfft} by means of the NFFT yields better results in this setting, the approximation~\eqref{eq:approx_exp_nfft_like} by means of our newly proposed method easily outperforms the classical NFFT in terms of the approximation error~\eqref{eq:maxerr_comparison_nfft}.

That is to say, Figure~\ref{fig:comparison_nfft} demonstrates that the novel NFFT-like approach in Algorithm~\ref{alg:nfft_generalized} is better suited for band\-limited functions,
while this superiority is not limited to~\mbox{$k\in\I_{M}$} but extends to the entire domain~\mbox{$v\in\big[-\frac M2,\frac M2\big]$}. Moreover, the error of Algorithm~\ref{alg:nfft_generalized} is bounded by the error estimates of the regularized Shannon sampling formulas in Section~\ref{sec:shannon}, whereas the quadrature error of the NFFT remains unclear.

\medskip
\begin{Example}
	\label{ex:nfft_generalized}
	To finally examine the approximation quality of the NFFT-like procedure in Algorithm~\ref{alg:nfft_generalized} for band\-limited functions we provide a function~$f$ with its corresponding Fourier transform~$\hat f$ in~\eqref{eq:inverse_integral}, such that we have access to the exact values~\mbox{$\hat f(k)$}, \mbox{$k\in\I_M$}, as input for Algorithm~\ref{alg:nfft_generalized}, as well as the exact function evaluations~\mbox{$f(x_j)$}, \mbox{$j=1,\dots,N$}.
	In doing so, we can compare the result~\mbox{$f_j$}, \mbox{$j=1,\dots,N$}, of Algorithm~\ref{alg:nfft_generalized} to the exact function evaluations~\mbox{$f(x_j)$}, \mbox{$j=1,\dots,N$}, by computing the maximum approximation error 
	\begin{align}
		\label{eq:err_NFFT_like}
		\max_{j=1,\dots,N} |f_j - f(x_j)| . 
	\end{align}
	For comparison we also compute the approximation error~\eqref{eq:err_NFFT_like} when~$f_j$ is the result of the classical NFFT in Algorithm~\ref{alg:nfft}.
	
	We consider the one-dimensional setting~\mbox{$d=1$} and for several band\-width parameters~\mbox{$M \in \{20,40,\dots,1000\}$} we study the function~\mbox{$f(x) = \sinc^2 \big(\frac{M}{2}\pi x\big)$} with the Fourier transform
	\begin{align*}
		\hat f(v)
		=
		\frac{2}{M} \cdot
		\begin{cases}
			1-\big|\frac{2v}{M}\big| & \colon |v| \leq \frac{M}{2}, \\
			0 & \colon \text{otherwise} .
		\end{cases}
	\end{align*}
	Note that the function~$f$ is scaled such that~\mbox{$\max_{x\in\R} f(x) = 1$} in\-de\-pend\-ent of the band\-width~$M$ and thereby the approximation errors~\eqref{eq:err_NFFT_like} are comparable for all considered~$M$.
	As evaluation points~\mbox{$x_j\in\big[-\frac 12+\frac mL,\, \frac 12-\frac mL\big]$}, \mbox{$j=1,\dots,N$}, we choose the scaled Chebyshev nodes
	\begin{align}
		\label{eq:points_cheb_scaled}
		x_j = \cos\left(\frac{(j-1)\pi}{N}\right) \cdot \left(\frac 12-\frac mL\right) ,
		\quad j=1,\dots,N ,
	\end{align}
	with~\mbox{$N = \frac{M}{2}$}, \mbox{$m=5$}, as well as~\mbox{$M_\sigma=L=M(1+\lambda)$} with~\mbox{$\lambda=1$}, and we use the $\sinh$-type window function~\eqref{eq:varphisinh}.
	
	The corresponding results are depicted in Figure~\ref{fig:approx_nfft_like}.
	As expected by Figure~\ref{fig:comparison_nfft}, the new NFFT-like procedure in Algorithm~\ref{alg:nfft_generalized} performs much better than the classical NFFT in Algorithm~\ref{alg:nfft}.
	While for~\mbox{$M\leq 80$} both approaches exhibit the same maximum approximation error~\eqref{eq:err_NFFT_like}, for larger band\-width~$M$ the approximation error~\eqref{eq:err_NFFT_like} gets smaller only for the NFFT-like procedure in Algorithm~\ref{alg:nfft_generalized}.
	That is to say, when approximating the evaluations~\mbox{$f(x_j)$}, \mbox{$j=1,\dots,N$}, of the band\-limited function~\mbox{$f\in\mathcal{B}_{M/2}(\R)$} by given samples~\mbox{$\hat f(k)$}, \mbox{$k\in\I_{M}$}, of the corresponding Fourier transform~\eqref{eq:inverse_integral}, reasonable results can be obtained by the NFFT in Algorithm~\ref{alg:nfft}, yet evidence indicates that our newly proposed NFFT-like procedure for band\-limited functions in Algorithm~\ref{alg:nfft_generalized} yields results that are at least as good, if not superior.
	Accordingly, we conclude that the NFFT-like procedure in Algorithm~\ref{alg:nfft_generalized} is the preferred approach in this context.
	\begin{figure}[ht]
		\centering
		\captionsetup[subfigure]{justification=centering}
		\begin{subfigure}[t]{0.45\columnwidth}
			\includegraphics[width=\textwidth]{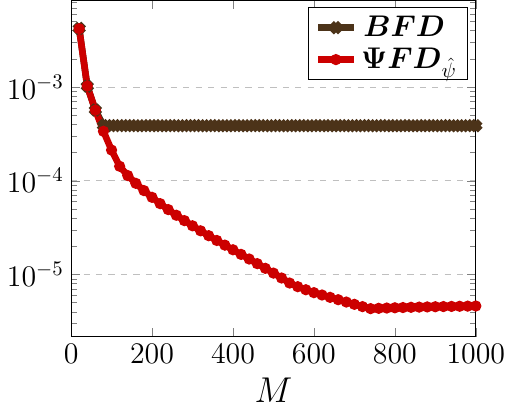}
		\end{subfigure}
		\caption{Maximum approximation error~\eqref{eq:err_NFFT_like} of Algorithms~\ref{alg:nfft} and~\ref{alg:nfft_generalized} using the $\sinh$-type window function~\eqref{eq:varphisinh} computed for the function~\mbox{$f(x) = \sinc^2 \big(\frac{M}{2}\pi x\big)$}, \mbox{$M \in \{20,40,\dots,1000\}$}, and the scaled Chebyshev nodes~\eqref{eq:points_cheb_scaled} with~\mbox{$N = \frac{M}{2}$}, \mbox{$m=5$}, \mbox{$M_\sigma=L=(1+\lambda)M$}, as well as~\mbox{$\lambda= 1$} and~\mbox{$d=1$}.
			\label{fig:approx_nfft_like}}
	\end{figure}
\end{Example}


\section*{Acknowledgments}
Melanie Kircheis acknowledges the support from the BMBF grant 01$\mid$S20053A (project SA$\ell$E) and the Deutsche Forschungsgemeinschaft (DFG, German Research Foundation) -- Project-ID 519323897.
Moreover, the authors thank the referee and the editor for their very useful suggestions for improvements.

\bibliographystyle{abbrv}

\end{document}